\documentclass{article}
\usepackage{amsmath,amssymb}
\usepackage{graphicx,epsfig}

\title{A note on extended full waveform inversion}
\author{Tristan van Leeuwen \\ Utrecht University, Utrecht, The Netherlands}

\begin{document}

\maketitle

\begin{abstract}
Full waveform inversion (FWI) aims at estimating subsurface medium properties from measured seismic data. It is usually cast as a non-linear least-squares problem that incorporates uncertainties in the measurements. In exploration seismology, extended formulations of FWI that allow for uncertaties in the physics have been proposed. Even when the physics is modelled accurately, these extensions have been shown to be beneficial because they reduce the non-lineary of the resulting data-fitting problem.
In this note, I derive an alternative (but equivalent) formulation of extended full waveform inversion. This re-formulation takes the form of a conventional FWI formulation that includes a medium-dependent weight on the residuals. I discuss the implications of this re-formulation and illustrate its properties with a simple numerical example.
\end{abstract}

\section{Introduction}
The seismic acquisition process can be described in terms of a process and measurement model:
\begin{eqnarray}
A(\mathbf{m})\mathbf{u} & = & \mathbf{q} + \boldsymbol{\eta},\\
P\mathbf{u}             & = & \mathbf{d} + \boldsymbol{\epsilon},
\end{eqnarray}
where $\mathbf{u}$ denotes the wavefield, $P$ is the sampling operator, $\mathbf{d}$ the observed data, $\mathbf{m}$ the medium parameters, $A$ the (discretized) wave-equation operator, $\mathbf{q}$ the source term and $\boldsymbol{\epsilon}, \boldsymbol{\eta}$ represent uncertainties in the models.

In full waveform inversion one tries to find a pair $(\mathbf{m},\mathbf{u})$ that fit the data and obeys the physics \cite{Tarantola1984}. Modelling the uncertainty in a Bayesian fashion, this leads to a Maximum a posteriori (MAP) estimation problem. In case $\boldsymbol{\epsilon}, \boldsymbol{\eta}$ are independently normally distributed with zero mean and covariances $\Sigma_m, \Sigma_p$ this leads to a non-linear least-squares problem
\begin{equation}
\label{eq:joint}
\min_{\mathbf{m},\mathbf{u}} \|P\mathbf{u} - \mathbf{d}\|_{\Sigma_{m}}^2 + \|A(\mathbf{m})\mathbf{u} - \mathbf{q}\|_{\Sigma_{p}}^2,
\end{equation}
where $\|\mathbf{r}\|_{W}^2 = \mathbf{r}^*W^{-1}\mathbf{r}$ denotes a weighted norm \cite{Tarantola2005}. When the uncertainty in the process model is neglegible (i.e., we have $\Sigma_p  = \sigma_p I$ with $\sigma_p \rightarrow 0$) we can eliminate $\mathbf{u} = A(\mathbf{m})^{-1}\mathbf{q}$ and retrieve the conventional formulation of full-waveform inversion:
$$\min_{\mathbf{m}} \|PA(\mathbf{m})^{-1}\mathbf{q} - \mathbf{d}\|_{\Sigma_m}.$$

Even if the process model is adequate, extended formulations of FWI of the form \eqref{eq:joint} have proven usefull in reducing the non-linearity of the problem \cite{vanLeeuwen2013Penalty1,Berkhout2014,van2015penalty,Huang2018,Huang2018a}. Some of these extensions, like the \emph{extened source} formulation \cite{Huang2018} or the \emph{contrast-source} formulation \cite{berg97} are closely related to \eqref{eq:joint} (see appendix \ref{sec:appext} for details).

In this note, I derive an alternative (but equivalent) reduced formulation of \eqref{eq:joint} of the form
\begin{equation}
\label{eq:mainresult}
\min_{\mathbf{m}} \|PA(\mathbf{m})^{-1}\mathbf{q} - \mathbf{d}\|_{\Sigma(\mathbf{m}) + \Sigma_m}^2,
\end{equation}
that includes a parameter-dependent covariance matrix that weighs the residual. In the Bayesian framework, we can iterpret $\Sigma(\mathbf{m}) + \Sigma_m$ as the covariance of the marginal of the posterior distribution. Indeed, \cite{Fang2018} propose to use this marginal for uncertainty quantification, but the exact form presented here seems to have been overlooked.

The alternative formulation \eqref{eq:mainresult} allows us to analyze the limit of vanishing measurement uncertainty; it immediately suggests that the residual should be measured in the norm $\|\cdot\|_{\Sigma(\mathbf{m})}^2$. While the limit of vanishing \emph{process} uncertainty has been investigated previously \cite{van2015penalty}, this result for vanishing \emph{measurement} uncertainty is novel.
Furthermore, this formulation may be more convenient for implementation since it only requires an additional weight matrix to be applied to the residual at each iteration of a conventional FWI workflow.

The remainder of the note is organized as follows. First, I give a derivation and disussion of the main result. An expression for the gradient of the new objective function is included. Then, I illustrate the behaviour of the new objective for a simple toy problem. Finally, I discuss some directions for future research and present conclusions.

\section{Theory}
In this section I sketch the main steps involved in deriving \eqref{eq:mainresult} from \eqref{eq:joint}.
The first step is to eliminate the state, $\mathbf{u}$, from \eqref{eq:joint} to obtain a reduced formulation. As the objective is quadratic in $\mathbf{u}$ we get a closed-form solution by solving the normal equations
\begin{eqnarray}
\label{eq:u}
\mathbf{u}(\mathbf{m}) = \left( A^*(\mathbf{m})\Sigma_p^{-1}A(\mathbf{m}) + P^*\Sigma_{m}^{-1}P\right)^{-1} \times  \nonumber\\\left(A^*(\mathbf{m})\Sigma_p^{-1}\mathbf{q} + P^*\Sigma_{m}^{-1}\mathbf{d}\right),
\end{eqnarray}
where $^*$ denotes the adjoint of an operator / complex-conjugate-transpose of a vector. Plugging this back in \eqref{eq:joint} and re-organizing terms (see appendix \ref{sec:appmain} for details) gives
\[
\phi(\mathbf{m}) = \|PA(\mathbf{m})^{-1} - \mathbf{d}\|_{\Sigma(\mathbf{m}) + \Sigma_m}^2,
\]
with
\[
\Sigma(\mathbf{m}) = P\left(A^*(\mathbf{m})\Sigma_p^{-1}A(\mathbf{m})\right)^{-1}P^*.
\]
This constitutes a parameter-dependent residual-weighting which can be thought of as the correlation of receiver-side Green's functions. Computing the full matrix would require a number of wave-simulations equal to the number of receivers. In practice, however, it may be possible to usefully approximate it or re-use some of the computations. Further discussion of these issues is postponed to a later section.

Remarkably, the gradient of the objective has a simple expression\footnote{This expression may be derived using the quotient-rule for matrix-differentiation: $\partial A^{-1} = - A^{-1}\left(\partial A\right)A^{-1}$.}, similar to that of conventional FWI:
\[
\frac{\partial\phi(\mathbf{m})}{\partial m_k} = -2\mathbf{u}_0^*\left(\frac{\partial A}{m_k}\right)\mathbf{v}_0 + 2 \mathbf{v}_0^*\left(\frac{\partial A}{m_k}\right)\mathbf{w}_{0},
\]
with $\mathbf{u}_0 = A^{-1}\mathbf{q}$, $\mathbf{v}_0 = A^{-*}P^*(\Sigma^{-1} + \Sigma_{m}^{-1})(P\mathbf{u}_0 - \mathbf{d})$ and $\mathbf{w}_0 = A^{-1}\Sigma_p^{-1}\mathbf{v}$. Note that we need only one additional forward solve to compute the gradient.

An interesting connection arises when consider the case of invertible $P$. In this case we find
\[
\phi(\mathbf{m}) = \|A(\mathbf{m})P^{-1}\mathbf{d} - \mathbf{q}\|_{\Sigma_p}^2.
\]
The depedency on the parameter is through $A$ in stead of $A^{-1}$, arguably making the problem easier to solve. This approach to estimating parameters from (nearly) full measurements of the state is sometimes referred to as the \emph{equation error method} \cite{Banerjee2013} and is closely related to wavefield gradiometry \cite{Poppeliers2013,DeRidder2017}. When $P$ is not invertible we cannot easily analyse the extended formulation in this fashion, but I will illustrate it with an example in the next section.

\section{Example}
\label{sec:example}
As an example, consider the constant-coefficient wave-equation in $\mathbb{R}^3$:
\[
(\partial_t^2 - c^2\nabla^2)u(t,\mathbf{x}) = q(t)\delta(\mathbf{x} - \mathbf{x}_s),
\]
with $u(t,\mathbf{x}) = 0$ for $t < 0$. The solution is given by
\[
u(t,\mathbf{x}) = Gq(t,\mathbf{x}) = \int g(t-t', \mathbf{x} - \mathbf{x}_s) q(t') \mathrm{d}t',
\]
with $g(t,\mathbf{x}) = \frac{\delta(t - \|\mathbf{x}\|_2 / c)}{\|\mathbf{x}\|_2}$. This leads to a time-extended formulation, which we can think of as the general extended formulation \eqref{eq:joint} where $\Sigma_p$ represents a delta-pulse in space.

The adjoint of the forward opertor is given by
\[
G^*u(t) = \int\!\!\int g(t'-t, \mathbf{x}' - \mathbf{x}_s) u(t',\mathbf{x}') \mathrm{d}t'\mathrm{d}\mathbf{x}'.
\]
The kernel $GG^*$ is then
%
%
involves convolution with
\[
k(t,t',\mathbf{x},\mathbf{x}') = \int g(t-t'', \mathbf{x} - \mathbf{x}_s) g(t'-t'', \mathbf{x}' - \mathbf{x}_s) \mathrm{d}t''.
\]
This yields
\[
k(t,t',\mathbf{x},\mathbf{x}') = \frac{\delta\left((t - t') - (\|\mathbf{x}-\mathbf{x}_s\|_2 - \|\mathbf{x}'-\mathbf{x}_s\|_2)/c\right)}{\|\mathbf{x}-\mathbf{x}_s\|_2\|\mathbf{x}'-\mathbf{x}_s\|_2}.
\]
For a single receiver-receiver pair, the convolution kernel $k(t,t',\mathbf{x}_{r,i}, \mathbf{x}_{r,j})$
thus contains one event with a time-lag corresonding to the difference in offset between the two receivers.


As a numerical example, we consider measuring the response of a single point-source at three receiver locations, at distances $0.8, 1$ and $1.2$ km from the source. The corresponding forward operator $F(c) PG(c)\mathbf{q}$ is implemented using a fast Fourier transform. The objective function reads
\[
\phi(c) = \|F(c)\mathbf{q}  - \mathbf{d}\|_{\sigma_m^2 I + \sigma_p^2 K(c)}^2,
\]
with $K(c) = F(c)F(c)^*$. Evaluating $\phi$ entails two steps; \emph{i)} compute the residual $\mathbf{r} = F(c)\mathbf{q}  - \mathbf{d}$, \emph{ii)} solve a (regularized) multidimensional deconvolution (MDD) problem $\left(\sigma_p^2 K(c) + \sigma_m^2 I\right)\widetilde{\mathbf{r}} = \mathbf{r}$. The objective value is obtained by computing the inner product $\phi(c) = \mathbf{r}^*\widetilde{\mathbf{r}}$. In this example, I solve the MDD problem using LSQR.

\begin{figure}
\includegraphics[scale=.5]{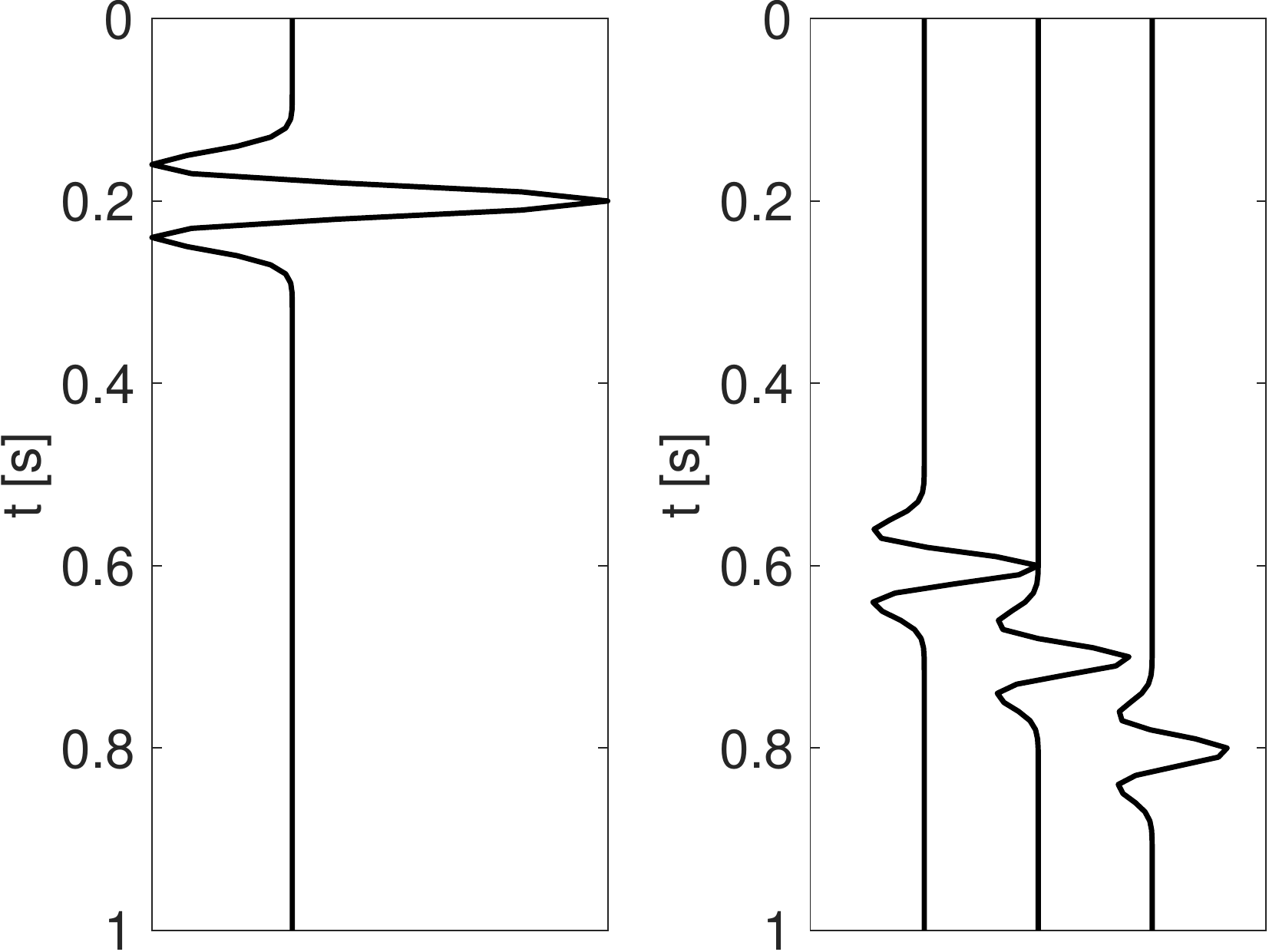}
\caption{Source and data for $c = 2$ km/s.}
\label{fig:example1a}
\end{figure}

\begin{figure}
\includegraphics[scale=.5]{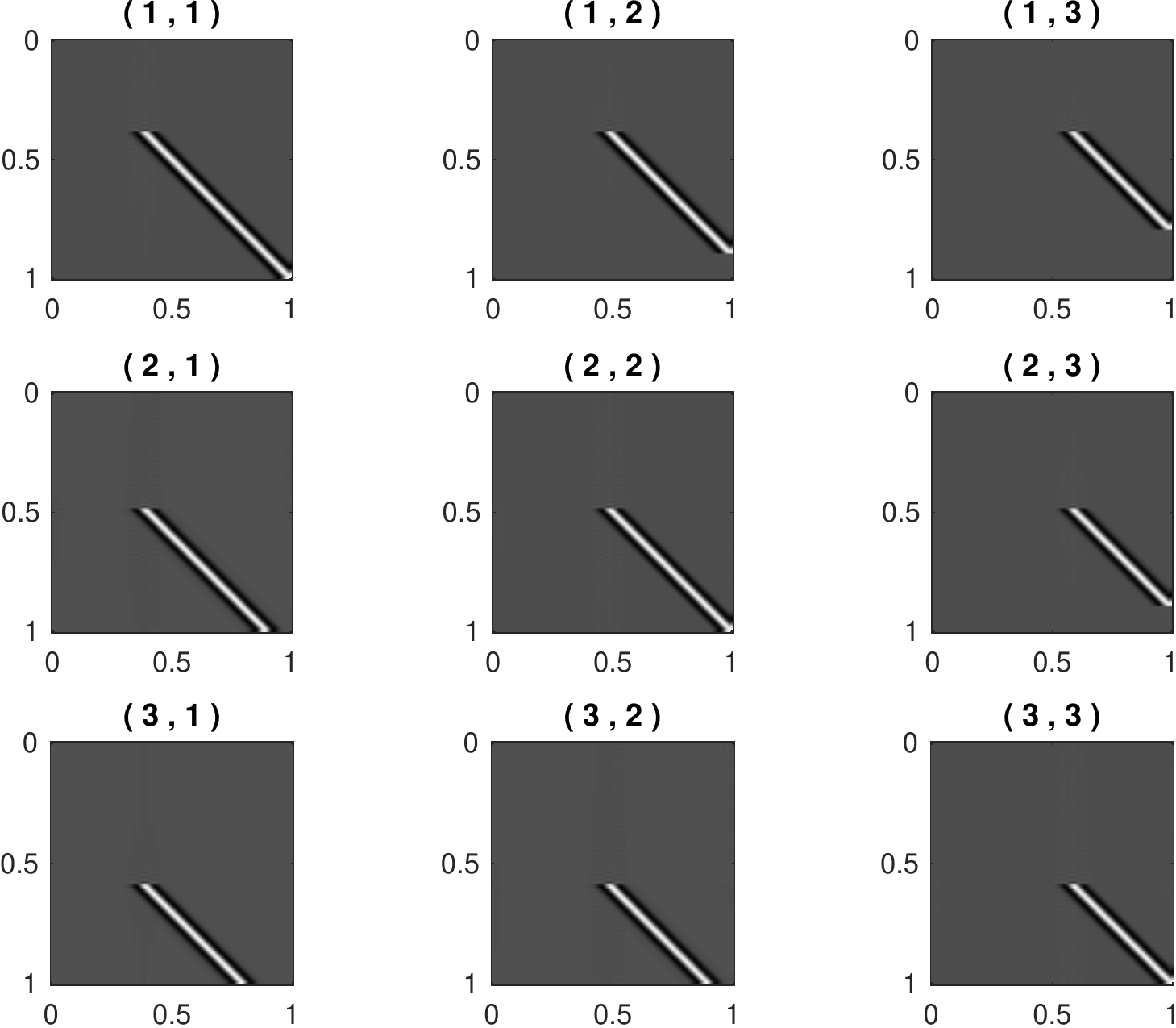}
\caption{Kernel $k(t,t',\mathbf{x}_r, \mathbf{x}_r')$ for $c = 2$ km/s. Depicted are nine panels corresponding to $k(t,t',\mathbf{x}_{r,i}, \mathbf{x}_{r,j})$ for $i,j \in \{1,2,3\}$.}
\label{fig:example1b}
\end{figure}

\begin{figure}
\includegraphics[scale=.25]{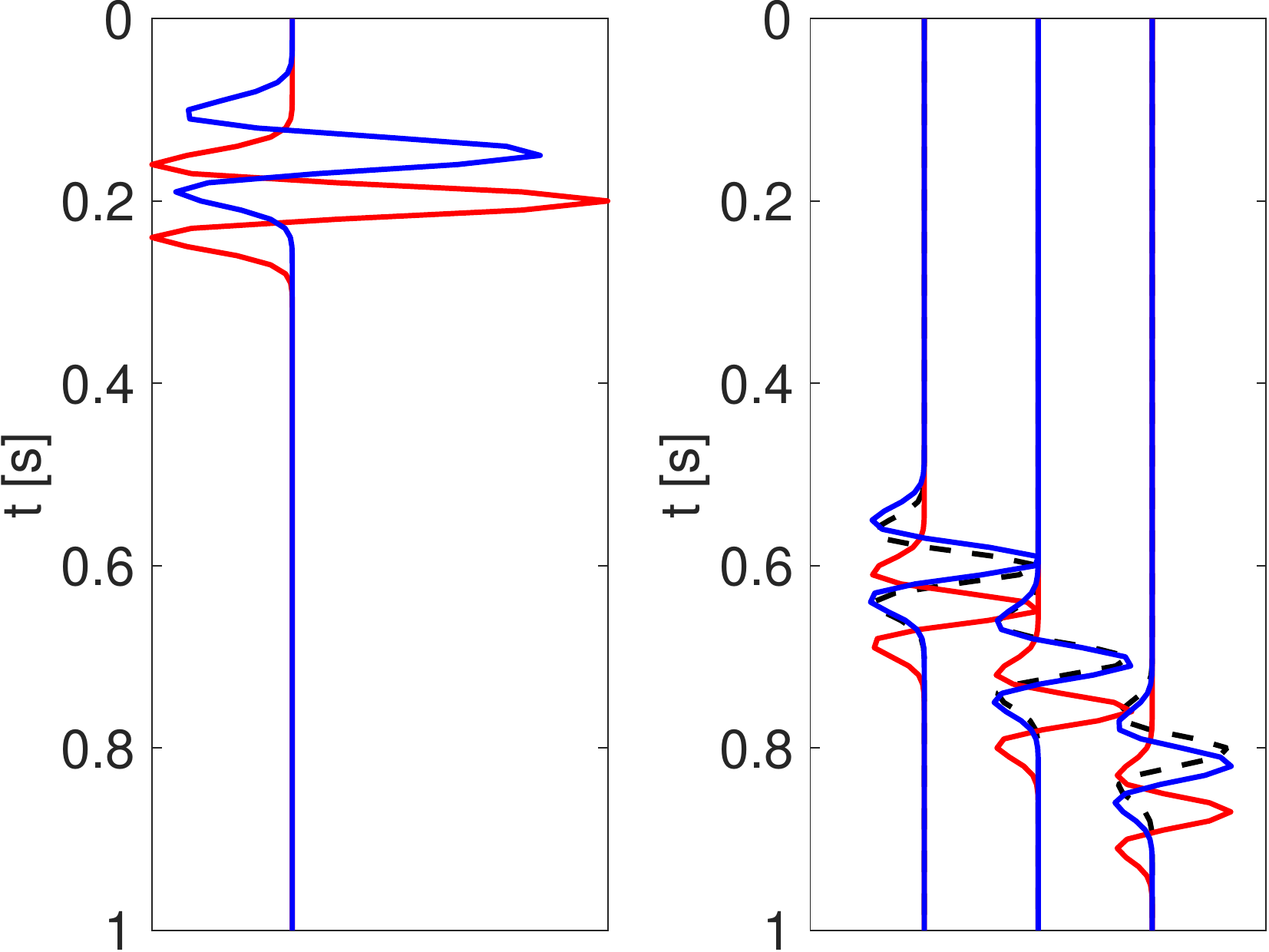}
\includegraphics[scale=.25]{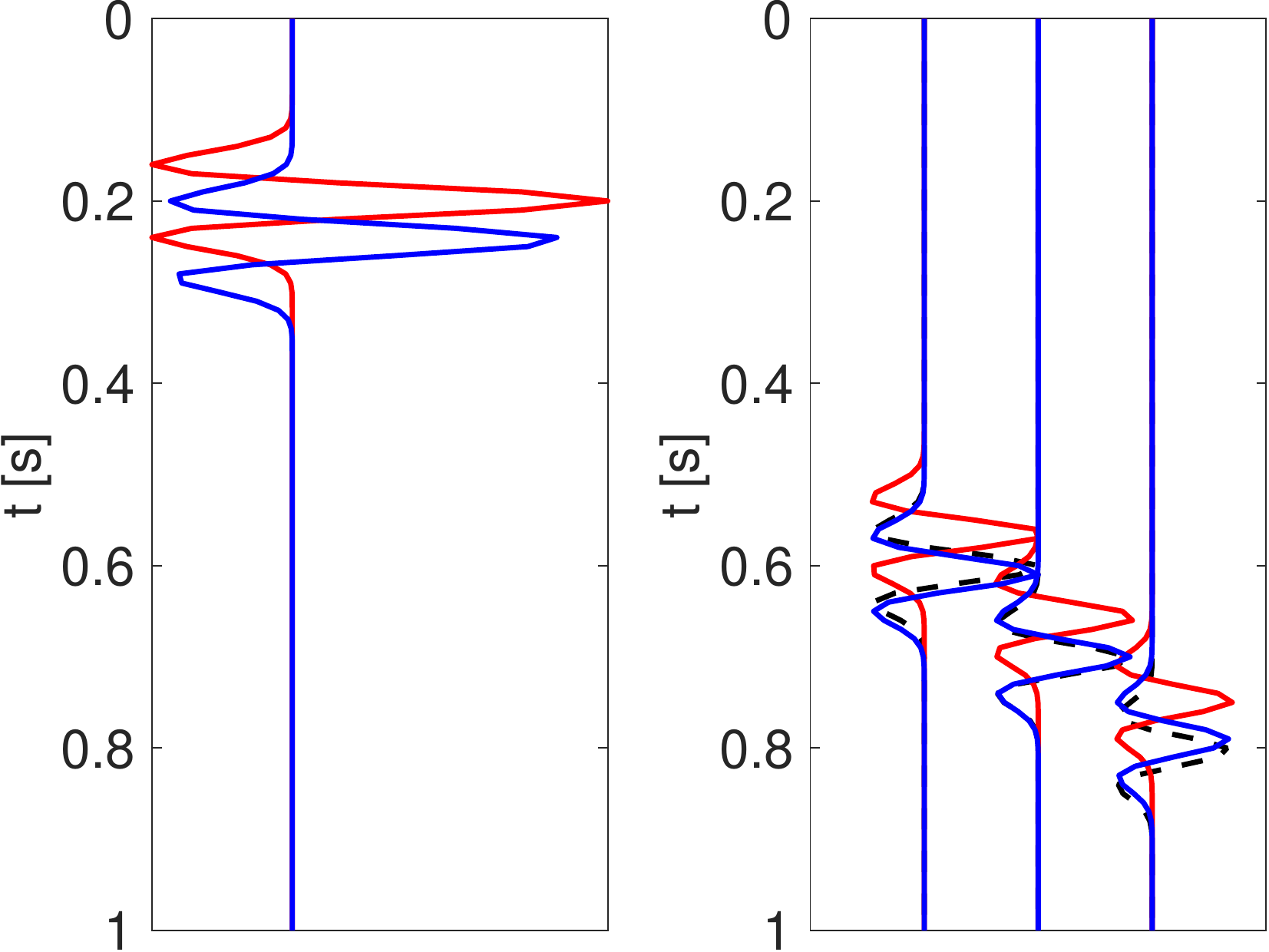}
\caption{Shown here are $\mathbf{q}$, $G(c)\mathbf{q}$ (red) and the estimated extended source $\mathbf{q} + \mathbf{r}$ and the corresponding data $G(c)(\mathbf{q} + \mathbf{r})$ (blue) for $c = 1.8$ km/s (left) and $c = 2.2$ km/s (right). We see that the extended source compensates for the velocity difference to fit the data (shown in black).}
\label{fig:example2}
\end{figure}

\begin{figure}
\includegraphics[scale=.3]{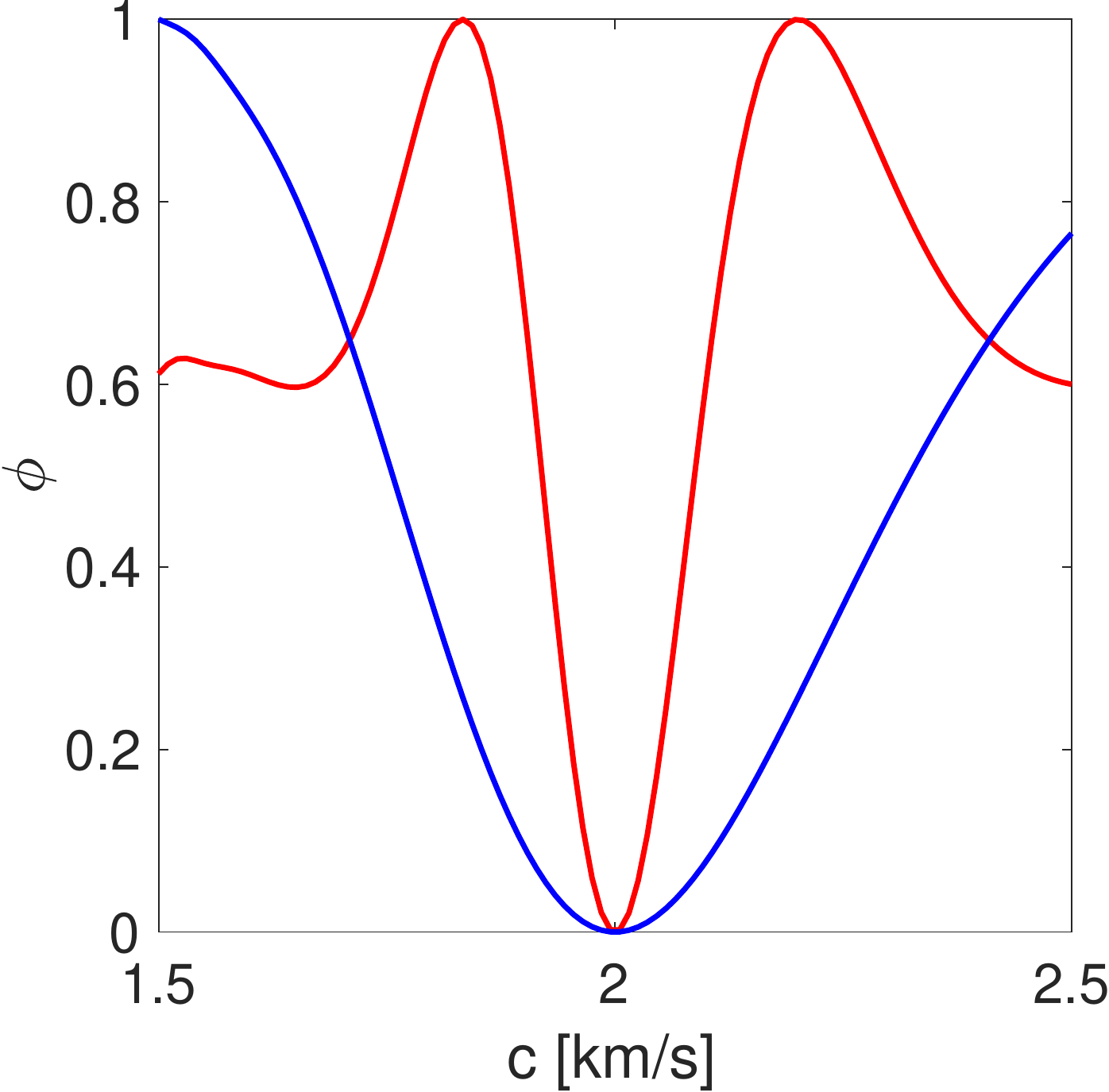}
\caption{Reduced objective, $\phi(c)$, for two limiting cases $\sigma_p\rightarrow 0$ (blue) and $\sigma_m\rightarrow 0$ (red)}
\label{fig:example3}
\end{figure}

The source function $\mathbf{q}$ and corresponding data for $c = 2$ km/s are depicted in figure \ref{fig:example1a}. This will serve as the ground truth solution that we aim to invert for. The corresponding kernel $K(c)$ is shown in figure \ref{fig:example1b}. We can clearly identity the events corresponding to the distance between the receivers. Figure \ref{fig:example2} shows the estimated extened sources for $c = 1.8$ km/s and $c = 2.2$ km/s and the corresponding data. We clearly see the difference between the two limiting cases; for $\sigma_p\rightarrow 0$ (conventional) we retreive the original source and do not fit the data, whereas for $\sigma_m \rightarrow 0$
(extended), additional events are being introduced in the source to fit the data. Finally,
the reduced objective $\phi(c)$ is depicted in figure \ref{fig:example3} in two limiting cases $\sigma_p\rightarrow 0$ (conventional) and $\sigma_m \rightarrow 0$ (extended). It clearly shows the reduced non-linearity of the extended formulation.

\section{Discussion}
The equivalence of the extended formulation of FWI and a residual-weighted version of conventional FWI opens up new possibilities for implementing the former. Note that the derivation presented here is valid for both time and frequency domain formulations of FWI. Of course, one would not explicitly form the system matrix $A$ for a time-domain formulation, but the action of the opertor, its adjoint and its inverse can be readily computed. Software libraries like RVL \cite{Padula2009}, SPOT \cite{SPOT}, ODL \cite{ODL}, PyLops \cite{PyLOps} and DeVito \cite{Louboutin2017} allow one to efficiently expose both time and frequency domain propagators as linear operators.

While forming and inverting the full model-dependent covariance matrix at each iteration may not be feasible in practice, various approximations may be investigated. Three flavours come to mind:
\begin{itemize}
\item Randomized matrix-probing techniques may be used to obtain a structured approximation of $\Sigma(\mathbf{m})$ that can be efficiently inverted. An extensive overview of such methods is given by \cite{DBLP:journals/corr/Wang15h}. It may not be necessary to construct such an approximation at each iteration, so that the cost of constructing it can be amortized over a number of iterations.
\item The data computed at the current iterate can be used to construct a covariance matrix: $\Sigma(\mathbf{m}) = \sum_i \mathbf{d}_i(\mathbf{m})\mathbf{d}_i(\mathbf{m}) = \sum_i PA(\mathbf{m})^{-1}\mathbf{q}_i\mathbf{q}_i^*A(\mathbf{m})^{-T}P^*$, leading to an effective process-model uncertainty given by $\Sigma_p = \left(\sum_i \mathbf{q}_i\mathbf{q}_i^*\right)^{-1}$. This can be interpreted as assinging large uncertainty to the region near the sources and no uncertainty away from it. Similarly, the measurement covariance matrix can be estimated by cross-correlating the residuals at each iterate \cite{Aravkin2012c}.
\item For simple velocity models, (semi-)analytic expressions of $\Sigma(\mathbf{m})$ may be derived, either using exact Green's functions for certain velocity models or using asymptotic (ray-based) approximations.
\end{itemize}

Of course, many FWI workflows contain data- or model-dependent weights on the residual (e.g, offset-weighting, amplitude balancing, etc.). The specific choice of using $\Sigma(\mathbf{m})$ as weight, however, is consistent with an extended formulation of FWI. This connectin may aid in further analysis of both conventional and fully extended FWI as both are now seen as extreme cases of the same formulation.

Note that the kernel expressed in section \ref{sec:example} suggests that it may be estimated from the data directly by correlating the traces corresponding to the same source. This would give us access to the kernel corresponding to the true model and information may be gleaned from that directly. This connection may shed new light on the extended formulation and warrants further investigation. It also points to a tentative connection to the work of \cite{Mamonov2015, Borcea2018} on model-order reduction for seismic imaging, where a similar kernel matrix is estimated to the data directly.


\section{Conclusion}
In this note I have derived a model-dependent residual-weight for full waveform inversion that makes it equivalent to an extended formulation. The residual and gradient can be readily computed using the standard tools available in any FWI workflow (i.e., forward and adjoint simulations, zero-lag correlations of wavefields). The effect of the residual weighting is shown on a simple toy-problem and appears to mitigate the non-linearity of the problem somewhat.


\appendix

\section{Extended formulations}
\label{sec:appext}
Starting from \eqref{eq:joint}
\[
\min_{\mathbf{m},\mathbf{u}}\|P\mathbf{u} - \mathbf{d}\|_{\Sigma_{m}}^2 + \|A\mathbf{u} - \mathbf{q}\|_{\Sigma_{p}}^2,
\]
I will point out connections to similar extended formulations of FWI.

The extended source formulation \cite{Huang2018} can be obtained by introducing a new variable $\mathbf{f} = A\mathbf{u} - \mathbf{q}$. This yields $\mathbf{u} = A^{-1}(\mathbf{m})(\mathbf{q} + \mathbf{f})$, and we get
\[
\min_{\mathbf{m},\mathbf{f}}\|PA^{-1}(\mathbf{m})(\mathbf{q} + \mathbf{f}) - \mathbf{d}\|_{\Sigma_{m}}^2 + \|\mathbf{f}\|_{\Sigma_{p}}^2.
\]
Note that this extension involves a fully general spatio-temporal source function. In practice, one may want to restrict the allowed sources. A source that extends only spatially can be thought of as corresponding to choosing the process uncertainty, $\Sigma_p$, to represent a delta-pulse in time. The annihilator used by \cite{Huang2018} to focus the source can then be interpreted as imposing a spatially decaying uncertainty.

The \emph{contrast source} \cite{berg97} formulation can be expressed as
\[
\min_{\mathbf{m}, \mathbf{w}} \|PA_0^{-1}(\mathbf{w} + \mathbf{q}) - \mathbf{d}\|_{\Sigma_m}^2 + \|\Delta A(\mathbf{m})A_0^{-1}(\mathbf{w} + \mathbf{q}) - \mathbf{w}\|_{\Sigma_p}^2,
\]
where $A_0$ is the opertor for a fixed background velocity, $\mathbf{w}$ is referred to as the \emph{contrast source} and $\Delta A(\mathbf{m})$ is the contrast w.r.t. the background velocity.


\section{Proof of main result}
\label{sec:appmain}
Start from
\[
\phi(\mathbf{m}) = \|P\mathbf{u}(\mathbf{m}) - \mathbf{d}\|_{\Sigma_{m}}^2 + \|A\mathbf{u}(\mathbf{m}) - \mathbf{q}\|_{\Sigma_{p}}^2,
\]
with $\mathbf{u}(\mathbf{m})$ as defined in \eqref{eq:u}. In the following, I will leave out the dependence on $\mathbf{m}$ for ease of notation. Now, introduce new variables
\[
\mathbf{v} = \mathbf{u} - A^{-1}\mathbf{q}, \, \mathbf{r} = \mathbf{d} - PA^{-1}\mathbf{q}.
\]
We get
\[
\phi(\mathbf{m}) = \|P\mathbf{v} - \mathbf{r}\|_{\Sigma_{m}}^2 + \|A\mathbf{v}\|_{\Sigma_{p}}^2,
\]
with $\mathbf{v}$ defined as
\begin{eqnarray*}
\mathbf{v} &=& \left(P^*\Sigma_{m}^{-1}P + A^*\Sigma_{p}^{-1}A\right)^{-1}P^*\Sigma_m^{-1}\mathbf{r}\\
&=& \left(A^*\Sigma_{p}^{-1}A\right)^{-1}P^*\left(P\left(A^*\Sigma_{p}^{-1}A\right)^{-1}P^* + \Sigma_m\right)^{-1}\mathbf{r}.
\end{eqnarray*}
The first expression follows directly form the definitions, while the second follows from the matrix identity (6.525) in \cite{Tarantola2005}. With this we can express
\[
P\mathbf{v} = K(K + \Sigma_m)^{-1}\mathbf{r},
\]
\[
A\mathbf{v} = \Sigma_{p}A^{-T}P^*(K + \Sigma_m)^{-1}\mathbf{r},
\]
with $K = P(A^*\Sigma_{p}^{-1}A)^{-1}P^*$.
This yields
\[
\|A\mathbf{v}\|_{\Sigma_p}^2 = \mathbf{r}^*(K + \Sigma_m)^{-1}K(K + \Sigma_m)^{-1}\mathbf{r},
\]
and
\[
\|P\mathbf{v} - \mathbf{r}\|_{\Sigma_m}^2 = \mathbf{r}^*\left(K(K + \Sigma_m)^{-1} - I\right)^* \Sigma_{m}^{-1} \left(K(K + \Sigma_m)^{-1} - I\right)\mathbf{r}
\]
Assembling all the terms and factoring out $(K + \Sigma_m)^{-1}$ we get
\[
\phi(\mathbf{m}) = \mathbf{r}^*(K + \Sigma_m)^{-1}\left(K + \Sigma_m\right)(K + \Sigma_m)^{-1}\mathbf{r}.
\]
This reduces to
\[
\phi(\mathbf{m}) = \mathbf{r}^*\left(K + \Sigma_m\right)^{-1}\mathbf{r},
\]
giving the desired result.

\end{document}